\documentclass[final]{siamltex}
\usepackage{graphicx}
\usepackage{verbatim}
\usepackage{epstopdf}
\usepackage{multirow}
\usepackage{subfig}
\usepackage{xcolor}

\def\headerframe#1#2#3{\vbox{\hrule height#2\hbox{\vrule width#2\hskip#1
\vbox{\vskip#1#3\vskip#1}\hskip#1\vrule width#2}\hrule height#2}}

\title{Stochastic Data Clustering 
         }

\author{Carl D. Meyer \thanks{Department of Mathematics and Institute of Advanced Analyitics, North Carolina State University, Raleigh, NC 27695, USA}
\and Charles D. Wessell \thanks{Department of Mathematics, Gettysburg College, Gettysburg, PA 17325, USA (\texttt{cwessell@gettysburg.edu})}}

\begin{document}

\maketitle

\begin{abstract}
In 1961 Herbert Simon and Albert Ando published the theory behind the long-term behavior of  a dynamical system that can be described by a nearly uncoupled matrix. Over the past fifty years this theory has been used in a variety of contexts, including queueing theory, brain organization, and ecology. In all these applications, the structure of the system is known and the point of interest is the various stages the system passes through on its way to some long-term equilibrium.

This paper looks at this problem from the other direction. That is, we develop a technique for using  the evolution of the system to tell us about its initial structure, and we use this technique to develop a new algorithm for data clustering.
\end{abstract}



\pagestyle{myheadings}
\thispagestyle{plain}
\markboth{C. D. MEYER AND C. D. WESSELL}{STOCHASTIC DATA CLUSTERING}

\section{Introduction} \label{sec:intro}
There is no shortage of data clustering algorithms. Indeed, many individual algorithms provide one or more parameters that can be set to a variety of values, effectively turning that single algorithm into many. Even if we restrict ourselves to a single algorithm with fixed starting parameters, we can still get varied results since methods like $k$-means and nonnegative matrix factorization (NMF) use random initializations that can lead to different final results. 

Rather than be frustrated by this seeming inconsistency of solutions, some clustering researchers have approached this problem with the goal of using all these clusterings to arrive at a single clustering solution that is superior to any individual solution.

The purpose of this article is to motivate and develop a new method for merging multiple clustering results using theory on the behavior of nearly uncoupled matrices developed by Nobel laureate Herbert Simon and his student Albert Ando.

When a collection of clustering methods is used, the collection is called an ensemble, and so this process is sometimes referred to as {\em ensemble clustering}. Others use the term {\em cluster aggregation} \cite{gionis2007ca}. Since the goal is for these varied methods to come to some agreement, it is also sometimes known as {\em consensus clustering}, which will be the term used throughout this paper.

The starting point for any clustering method is an $m$-dimensional data set of $n$ elements. The data set can thus be stored as an $m \times n$ matrix $A$ where each column represents an element of the data set and each row contains the value of a particular attribute for each of the elements. If the assignment of clusters from a single run of a clustering algorithm is denoted by $\mathcal{C}_{k}$, then the input to any consensus method will be $\mathcal{C} = \{\mathcal{C}_{1}, \mathcal{C}_{2}, \dots , \mathcal{C}_{r} \}$.

One approach for solving this problem is attempting to find a clustering $\mathcal{C}^{*}$ that is as close as possible to all the $\mathcal{C}_{k}$'s. This is an optimization problem known as \textit{median partition}, and is known to be NP-complete.  A number of heuristics for the median partition problem exist. Discussion of these heuristics with comparisons and results on real-world data sets can be found in \cite{filkov, filkovskiena, goderfilkov}.

Other researchers have brought statistical techniques to bear on this problem, using bootstrapping or other more general resampling techniques to cluster subsets of the original data set, and then examining the results using some measure of consistency to settle on the final clustering \cite{fredjain2003, consensusmit}. 

Additional approaches include a consensus framework built on a variational Bayes mixture of Gaussians model \cite{danish} and using algorithms originally intended for rank aggregation problems \cite{ailon2008aii}.

Other approaches to this problem begin by storing the information from each $\mathcal{C}_{k}$ in an $n \times n$ adjacency matrix $A^{(k)}$ such that if data set elements $i$ and $j$ are in the same cluster according to $\mathcal{C}_{k}$, then $a^{(k)}_{ij}=1$, and $a^{(k)}_{ij}=0$ if they are not (in this paper we will define $a^{(k)}_{ii}=1$ for $i=1,2, \dots, n$). The collection of these $r$ adjacency matrices can be used to define a hypergraph which can then be partitioned (i.e. clustered) using known hypergraph partitioning algorithms \cite{strehl}.

Alternatively, this collection of adjacency matrices can be summed to form the consensus matrix $S$. Each entry $s_{ij}$ of $S$ denotes how many times elements $i$ and $j$ clustered together. For those who would prefer that all entries of $S$ lie in the interval $[0,1]$, $S$ can be defined as the sum of the adjacency matrices times $\frac{1}{r}$, resulting in a symmetric similarity matrix whose similarity measure is the fraction of the time that two elements were clustered together. In this paper, $S$ will always be used to refer to the sum of the adjacency matrices.

Once $S$ is constructed, its columns can be clustered and thus the original data is clustered \cite{race}. This method using single-link hierarchical clustering on $S$, after elements below a threshold have been zeroed out, has proven effective \cite{fredjain2002}.

A new methodology developed to cluster different conformations of a single drug molecule comes the closest to the approach developed in this paper. For this application, a Markov chain transition matrix can be created where the $ij$-th entry gives the probability the molecule changes from conformation $i$ to conformation $j$. The goal is to then find sets of conformations such that if the molecule is currently in a particular set, it will remain in that set for a relatively long time. Approaches to this clustering problem have included examination of the first few eigenvectors of the transition matrix (\cite{deuflhard2000iai} and then improved in \cite{deuflhard2005rbc}), clustering the data based on the second singular vector \cite{fritzsche2008sai,tifenbach2011sba}, and spectral analysis of a family of Hermitian matrices that is a function of the the transition matrix \cite{jacobi2010rsm}.  

\section{A new approach}
The data clustering method introduced in this paper is based on the 1950's variable aggregation work of the Nobel prize winning economist Herbert Simon and his graduate student Albert Ando \cite{simonando}. Their theory will be reviewed in Section \ref{sec:tb}, and further theoretical work will be developed in Sections \ref{sec:sk} -- \ref{sec:spec} before the algorithm is introduced in Section \ref{sec:cip}.

\subsection{Theoretical background} \label{sec:tb}
\begin{figure}[h] 
\begin{center}
\includegraphics[scale=.6,trim=0 150 0 0]{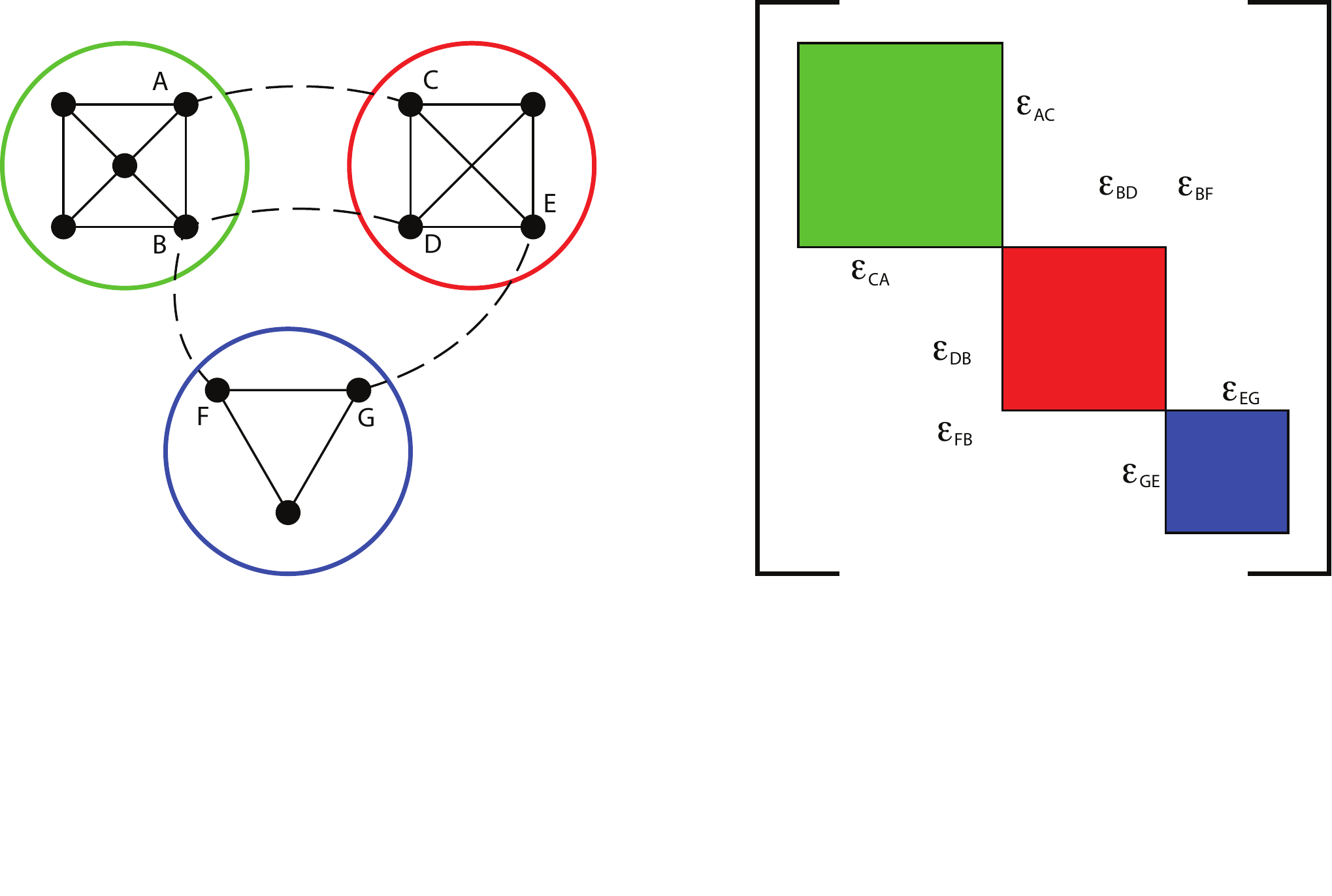}
\caption{This figure illustrates a simple Simon-Ando system and how it would be represented in matrix form. Let the circles on the left represent three small countries. The graphs within each circle represent companies in those countries and the solid lines between them represent a large amount of capital exchange between the companies. The dashed lines represent a small amount of cross-border exchange. A matrix whose entries represented the amount of economic activity between any two companies in this system would look like the one on the right with the shaded areas being dense with relatively large values and the epsilons being relatively small.} 
\label{fig:abc}
\end{center}   
\end{figure}

Simon-Ando theory was originally designed as a way of understanding the short and long term behavior of an economy with a certain structure. Figure \ref{fig:abc} illustrates a simple system where Simon-Ando theory would apply. 

Such a closed economic system, without any outside influences, is known to eventually reach a state of equilibrium, that is, after some initial fluctuations, the flow of goods and capital between any two industries will remain more or less constant. Rather than waiting for this economic equilibrium to occur, Simon and Ando tried to predict the long-term equilibrium by making only short-term observations. They proved that what happens in the short run completely determines the long-term equilibrium.

Over the years scholars in a variety of disciplines have realized the usefulness of a framework that represents a number of tightly-knit groups that have some loose association with each other, and Simon-Ando theory has been applied in areas as diverse as ecology \cite{levin1992pps}, computer queueing systems \cite{courtois}, brain organization \cite{sporns2000ccr}, and urban design \cite{salingaros2000cuc}. Simon himself went on to apply the theory to the evolution of multicellular organisms \cite{simon2002nds}.


The $n \times n$ matrix $S$ is called \emph{uncoupled} if it has the form
 \[ S=
\left(
\begin{array}{cccc}
S_{11}  & 0   & \dots &0   \\
0  & S_{22}   & \dots & 0    \\
\vdots   & \vdots   & \ddots & \vdots \\
0  & 0   & \dots & S_{kk}    \\ 
\end{array}
\right),
\]
where the diagonal blocks $S_{ii}$ are square. If $S$ is not uncoupled for any value of $k \ge 2$ and if entries in the off-diagonal blocks are small relative to those in the diagonal blocks, then we say that $S$ is \emph{nearly uncoupled}. The matrix in Figure \ref{fig:abc} is an example of a nearly uncoupled matrix. A more formal measure of uncoupledness will be introduced in Definition \ref{def:sigma}.

If the consensus matrix $S$ described in the Introduction is nearly uncoupled, we will show that Simon-Ando theory can be used to cluster the data it describes. Notice that $S$ is symmetric and this combined with it not being uncoupled means $S$ is also irreducible. For reasons that will soon become apparent, the new clustering method will require that $S$ be converted to doubly stochastic form. This new matrix will be called $P$ and the data clustering method will depend on $P$ having a unique stationary distribution vector (which is guaranteed by irreducibility) and a known structure (which is guaranteed by double stochasticity).

Before we can use $P$ to cluster data we need to introduce the concept of stochastic complementation.

If $P$ is stochastic then each  diagonal block $P_{ii}$ has a stochastic complement defined by
\begin{equation} \label{scdef}
C_{ii}=P_{ii}+P_{i\star} \left(I-P_{i}\right)^{-1}P_{\star i},
\end{equation}
where $P_{i}$ is the matrix obtained by deleting the $i$th row and $i$th column of blocks from $P$, $P_{i \star}$ is the $i$th row of blocks of $P$ with $P_{ii}$ removed, and $P_{\star i}$ is the $i$th column of blocks of $P$ with $P_{ii}$ removed. Since every principal submatrix of $I-P$ of order $n-1$ or smaller is a nonsingular $M$-matrix, the matrix $(I-P_{i})^{-1}$ found in (\ref{scdef}) is defined and $(I-P_{i})^{-1} \ge 0$. Furthermore, if $P$ is stochastic and irreducible, then each $C_{ii}$ is itself a stochastic, irreducible matrix with stationary distribution vector $c_{i}^{T}$ \cite{berplem, meyer}.

Let $x_{0}^{T}$ be a probability row vector and consider the evolution equation
\begin{equation} \label{1}
x_{t}^{T} = x_{t-1}^{T}P 
\end{equation}
or its equivalent formulation 
\begin{equation} \label{2}
x_{t}^{T} = x_{0}^{T}P^{t}. 
\end{equation}

Simon-Ando theory asserts that $x_{t}^{T}$ will pass through distinct stages as $t$ grows to infinity. Meyer \cite{meyer} describes how these stages can be interpreted in terms of the individual stationary distribution vectors $c_{i}^{T}$. The following lemma and theorem will aid in extending that explanation to the case where $P$ is doubly stochastic. The proof of the lemma is a direct application of principles of permutation matrices and is omitted.

\begin{lemma} \label{lemmaone} Let $P$ be an $n \times n$ irreducible doubly stochastic matrix in which the diagonal blocks are square. Let $Q$ be the permutation matrix associated with an interchange of the first and $i$th block rows (or block columns) and let $\tilde{P}$ be defined as
$$ \tilde{P} = QPQ.$$
If $\tilde{P}$ is partitioned into a $2 \times 2$ block matrix
\begin{equation} \label{twobytwo}
\tilde{P} = \left(
\begin{array}{cc}
\tilde{P}_{11}  & \tilde{P}_{12}    \\
\tilde{P}_{21}  & \tilde{P}_{22}  \\
\end{array}
\right)
\; where \; \tilde{P}_{11}=P_{ii},
\end{equation}
then the stochastic complement of $P_{ii}$ is
\begin{equation}
C_{ii}=\tilde{C}_{11}=\tilde{P}_{11}+\tilde{P}_{12} \left(I- \tilde{P}_{22}\right)^{-1}\tilde{P}_{21}
\end{equation}
\end{lemma}

\begin{theorem} If
 \[ P=
\left(
\begin{array}{cccc}
P_{11}  & P_{12}   & \dots & P_{1k}    \\
P_{21}  & P_{22}   & \dots & P_{2k}    \\
\vdots   & \vdots   & \ddots & \vdots \\
P_{k1}  & {P}_{k2}   & \dots & P_{kk}    \\ 
\end{array}
\right)
\]
is an irreducible doubly stochastic matrix, then each stochastic complement is also an irreducible, doubly stochastic matrix.
\end{theorem}

\begin{proof}
As stated earlier, if the stochastic matrix $P$ is irreducible, then so are each of its stochastic complements. Therefore, we need only prove that each $S_{ii}$ is doubly stochastic. For a given $i$, suppose diagonal block $P_{ii}$ has been repositioned such that $\tilde{P}_{11}=P_{ii}$ as in  (\ref{twobytwo}) of Lemma \ref{lemmaone}.

Let $e$ represent a column vector of all ones. Both the row and column sums of $P$ are one, so allowing the size of $e$ to be whatever is appropriate for the context, the following four equations are true 
\begin{equation}
\tilde{P}_{11}e + \tilde{P}_{12}e = e 
\end{equation}
\begin{equation} \label{this}
\tilde{P}_{21}e + \tilde{P}_{22}e = e 
\end{equation}
\begin{equation}
e^{T}\tilde{P}_{11} + e^{T}\tilde{P}_{21} = e^{T} 
\end{equation}
\begin{equation} \label{that}
e^{T}\tilde{P}_{12} + e^{T}\tilde{P}_{22} = e^{T}
\end{equation}

Equations \ref{this} and \ref{that} can be rewritten to yield
$$
e = \left(I - \tilde{P}_{22} \right) ^{-1} \tilde{P}_{21}e \;\; \mbox{and} \;\; e^{T}=e^{T} \tilde{P}_{12} \left(I - \tilde{P}_{22} \right)^{-1}.
$$
As noted earlier, $(I - \tilde{P}_{22}) ^{-1} \ge 0$, and hence 
$$
\tilde{C}_{11}= \tilde{P}_{11}+\tilde{P}_{12} \left(I- \tilde{P}_{22}\right)^{-1}\tilde{P}_{21} \ge 0.
$$
Multiplying $\tilde{C}_{11}$ on the right by $e$ and on the left by $e^{T}$ yields
$$
\tilde{C}_{11}e= \tilde{P}_{11}e+\tilde{P}_{12} \left(I- \tilde{P}_{22}\right)^{-1}\tilde{P}_{21}e =
 \tilde{P}_{11}e+\tilde{P}_{12}e = e
$$ 
and
$$ 
e^{T}\tilde{C}_{11}= e^{T}\tilde{P}_{11}+e^{T}\tilde{P}_{12} \left(I- \tilde{P}_{22}\right)^{-1}\tilde{P}_{21} =
e^{T}\tilde{P}_{11}+e^{T} \tilde{P}_{21} = e^{T}.
$$
Therefore, since $C_{ii}=\tilde{C}_{11}$, each stochastic complement is doubly stochastic.
\end{proof}

Markov chain theory tells us that as $t \rightarrow \infty$, $x_{t}^{T}$ will approach the uniform distribution vector $(1/n \;\; 1/n \;\; \dots \;\; 1/n)$. If the size of each $P_{ii}$ is $n_{i} \times n_{i}$, we also know that $c_{i}^{T} =  (1/n_{i} \;\; 1/n_{i} \;\; \dots \;\; 1/n_{i})$.

As $t$ increases from zero, $x_{t}^{T}$ initially goes through changes driven by the comparatively  large values in each $P_{ii}$. Once these changes have run their course, the system settles into a period of short-term stabilization characterized by
\begin{eqnarray*} \label{sts}
x_{t}^{T}  &\approx& ( \alpha_{1}c_{1} \;\;  \alpha_{2}c_{2} \;\; \dots \;\; \alpha_{k}c_{k}) \\
&=& \left(  
\frac{\alpha_{1}}{n_{1}} \, \frac{\alpha_{1}}{n_{1}} \, \dots \, \frac{\alpha_{1}}{n_{1}} \, \vline \,\frac{\alpha_{2}}{n_{2}} \, \frac{\alpha_{2}}{n_{2}} \, \dots \,  \frac{\alpha_{2}}{n_{2}} \, \vline \;\; \dots \;\; \vline \frac{\alpha_{k}}{n_{k}} \, \frac{\alpha_{k}}{n_{k}} \, \dots \, \frac{\alpha_{k}}{n_{k}}  \right) 
\end{eqnarray*}
where each $\alpha_{i}$ is a constant dependent on $x_{0}^{T}$.

After this equilibrium period, the elements of $x_{t}^{T}$ begin to change again through a period called middle-run evolution, this time being affected by the small values in the off-diagonal blocks, but the change is predictable and can be described by
\begin{eqnarray*} \label{mre}
x_{t}^{T}  &\approx& ( \beta_{1}c_{1} \;\;  \beta_{2}c_{2} \;\; \dots \;\; \beta_{k}c_{k}) \\
&=& \left( \frac{\beta_{1}}{n_{1}} \, \frac{\beta_{1}}{n_{1}} \, \dots \, \frac{\beta_{1}}{n_{1}} \,  \vline \,\frac{\beta_{2}}{n_{2}} \, \frac{\beta_{2}}{n_{2}} \, \dots \, \frac{\beta_{2}}{n_{2}} \, \vline \;\; \dots \;\; \vline \frac{\beta_{k}}{n_{k}} \, \frac{\beta_{k}}{n_{k}} \, \dots \, \frac{\beta_{k}}{n_{k}} \right)
\end{eqnarray*} 
where each $\beta_{i}$ is dependent on $t$.

Simon and Ando were not interested in clustering data. For them, the importance of stages like short-term stabilization and middle-run evolution lie in the fact that even for small values of $t$, the structure of ${x}_{t}^{T}$ reflected the stationary probability vectors of the smaller ${C}_{ii}$ matrices. From there, examination of the ${x}_{t}^{T}$ vector during the relatively stable periods would allow for determination of these smaller stationary probability vectors and facilitate the calculation of the stationary probability vector for $P$.

For cluster analysis however, the focus is turned around. Since we will be using doubly stochastic $P$ matrices, we already know that the stationary probability vector is the uniform probability vector. We also know that each diagonal block ${P}_{ii}$ is associated with a uniform probability vector related to its stochastic complement.  Identification of the clusters then comes down to examining the entries of ${x}_{t}^{T}$. The key is to look for elements of $x_{t}^{T}$ that are approximately equal. The only difference between short-run and middle-run is whether the elements of $x_{t}^{T}$ stay at approximately the same value for a number of iterations or move together towards the uniform probability distribution. 

All the development in this section assumed a doubly stochastic matrix. We will now consider how to convert a matrix into doubly stochastic form, and show that the process does not destroy any of the desirable characteristics of our matrix.

\subsection{Sinkhorn-Knopp} \label{sec:sk}
The process of converting a matrix into doubly stochastic form has drawn considerable attention, and in 1964 Sinkhorn showed that any positive square matrix can be scaled to a unique doubly stochastic matrix \cite{sinkhorn1964rap}. This result can be extended to nonnegative matrices as long as the zero entries are in just the right places. An understanding of this zero structure will require some definitions.

\begin{definition}
(Sinkhorn and Knopp \cite{sk}) A nonnegative $n \times n$ matrix $S$ is said to have \emph{total support} if $S \ne 0$ and if every  positive element of $S$ lies on a positive diagonal, where a \emph{diagonal} is defined as a sequence of elements $s_{1 \sigma(1)}, s_{2 \sigma(2)}, \dots, s_{n \sigma(n)}$ where $\sigma$ is a permutation of $\{1,2, \dots, n\}$.\footnote{Notice that by this definition of \emph{diagonal}, the main diagonal of a matrix is the one associated with the permutation $\sigma = ( 1 \; 2 \; 3\; \dots \; n)$.} 
\end{definition}

\begin{definition}
(Minc \cite{minc1988nm}, p.82) An $n \times n$ matrix $S$ is \emph{partly indecomposable} if there exist permutation matrices $P$ and $Q$ such that
\[
PSQ = 
\left[
\begin{array}{cc}
 X & Z     \\
 0 &   Y   
\end{array}
\right],
\]
where $X$ and $Y$ square. If no such $P$ and $Q$ exist, then $S$ is \emph{fully indecomposable}. 
\end{definition}

\begin{definition}
(Minc \cite{minc1988nm}, p.82) Two matrices $A$ and $B$ are \emph{permutation equivalent}, or \emph{p-equivalent}, if there exist permutation matrices $Q$ and $\hat{Q}$ such that $A=QB\hat{Q}$.
\end{definition}

This new terminology will help in understanding the following, nearly identical theorems that were independently proven and then published within a year of each other, the first in 1966 and the second in 1967.

\begin{theorem} \label{th:bps}
(Brualdi, Parter, and Schneider \cite{brualdi1966den}) If the $n \times n$ matrix $A$ is nonnegative and fully indecomposable, then there exist diagonal matrices $D_{1}$ and $D_{2}$ with positive diagonal entries such that $D_{1}AD_{2}$ is doubly stochastic. Moreover $D_{1}$ and $D_{2}$ are uniquely determined up to scalar multiples.
\end{theorem}

\begin{theorem} \label{th:sk}
(Sinkhorn and Knopp \cite{sk}) If the $n \times n$ matrix $A$ is nonnegative, then a necessary and sufficient condition that there exists a doubly stochastic matrix of the form $D_{1}AD_{2}$ where $D_{1}$ and $D_{2}$ are diagonal matrices with positive diagonal entries is that $A$ has total support. If $D_{1}AD_{2}$ exists, then it is unique. Also $D_{1}$ and $D_{2}$ are unique up to a scalar multiple if and only if $A$ is fully indecomposable.
\end{theorem}

The uniqueness up to a scalar multiple of $D_{1}$ and $D_{2}$ mentioned in both theorems means that if $E_{1}$ and $E_{2}$ are also diagonal matrices such that  $E_{1}AE_{2}$ is doubly stochastic, then $E_{1}=\alpha D_{1}$ and  $E_{2}=\beta D_{2}$ where $\alpha\beta=1$.

The way that the consensus similarity matrix $S$ is constructed guarantees its nonnegativity, so the only thing standing in the way of knowing that the scaling matrices $D_{1}$ and $D_{2}$ exist is showing that $S$ either has total support or is fully indecomposable. Reviewing the definitions of these terms, neither of these tasks seems inviting. Fortunately, there is a theorem that will simplify the matter.

\begin{theorem} \label{Minc}
(Minc \cite{minc1988nm}, p.86) A nonnegative matrix is fully indecomposable if and only if it is $p$-equivalent to an irreducible matrix with a positive main diagonal.
\end{theorem}

$S$ is trivially $p$-equivalent to itself since $S=ISI$ and $S$ is an irreducible matrix with a positive main diagonal. Now that we know $S$ is fully indecomposable, its symmetry is going to guarantee another useful result. The proof of the following lemma is included since there was a typographical error in the original paper. 

\begin{lemma} \label{Csima}
(Csima and Datta \cite{csima1972dts}) Let $S$ be a fully indecomposable symmetric matrix. Then there exists a diagonal matrix $D$ such that $DSD$ is doubly stochastic.
\end{lemma}

\begin{proof}
Let $D_{1}$ and $D_{2}$ be nonnegative diagonal matrices such that $D_{1}SD_{2}$ is doubly stochastic. Then $(D_{1}SD_{2})^{T} =  D_{2}SD_{1}$ is also doubly stochastic. By the uniqueness up to a scalar multiple from Theorems \ref{th:bps} and \ref{th:sk}, we know $D_{2}=\alpha D_{1}$ and $D_{1}=\beta D_{2}$. Using the first of these facts
\begin{eqnarray*}
D_{1}SD_{2} &=& D_{1}S\alpha D_{1} \\
&=& \sqrt{\alpha}D_{1}S\sqrt{\alpha}D_{1} \\
&=& DSD
\end{eqnarray*}
shows us that $D = \sqrt{\alpha}\;D_{1}$.
\end{proof}

\subsection{The structure of  $DSD$} \label{sec:dsd}

We will use $P$ as the symbol for the doubly stochastic matrix derived from $S$, that is $P=DSD$. For simplicity of notation, the $i^{th}$ diagonal entry of $D$ will be denoted $d_{i}$. We will show that $P$ has the same desirable properties that $S$ has. 

\begin{lemma} 
If $S$ is an $n \times n$ fully indecomposable irreducible matrix and $P=DSD$ is doubly stochastic, then $P$ is irreducible.
\end{lemma}

\begin{proof}
Since $S$ is irreducible, there is no permutation matrix $Q$ such that
\[
QSQ^{T} = 
\left[
\begin{array}{cc}
 X & Z     \\
 0 &   Y   
\end{array}
\right].
\]
where both $X$ and $Y$ are square.

Thus the only way that $P=DSD$ could be reducible is if the zero structure of $S$ is changed by the multiplication. But notice that since $p_{ij} = d_{i}d_{j}s_{ij}$ and both $d_{i}$ and $d_{j}$ are positive, $p_{ij}=0$ only when $s_{ij}=0$. So the zero structure does not change, and $P$ is irreducible. 
\end{proof}

Since the number of times elements $i$ and $j$ cluster with one another is necessarily equal to the number of times elements $j$ and $i$ cluster with one another, the symmetry of the consensus similarity matrix $S$ reflects a real-world property of the consensus clustering problem and so it is important that symmetry is not lost when $S$ is converted into $P$. 

\begin{lemma} 
If $S$ is an $n \times n$ fully indecomposable symmetric matrix and $P=DSD$ is doubly stochastic, then $P$ is symmetric.
\end{lemma}

\begin{proof}
\begin{equation}
P^{T} = (DSD)^{T} = DS^{T} D =  DSD = P 
\end{equation} 
\end{proof}

We wish to prove that if $S$ is nearly uncoupled, then so is $P$. To do so we first need a formal definition of near uncoupledness. Then we will show how this uncoupling measure for $P$ is related to the uncoupling measure of $S$.

\begin{definition} \label{def:sigma} 
Let $n_{1}$ and $n_{2}$ be fixed positive integers such that $n_{1}+n_{2}=n$, and let $S$ be an $n \times n$ symmetric, irreducible matrix whose respective rows and columns have been rearranged to the form
\[
S= \left[
\begin{array}{cc}
 S_{11} & S_{12}     \\
 S_{21} & S_{22}  
\end{array}
\right]
\]
where $S_{11}$ is $n_{1} \times n_{1}$  and $S_{22}$  is $n_{2} \times n_{2}$ so that the ratio
\[
\sigma(S,n_{1})=
\frac{e^{T}S_{12}e+e^{T}S_{21}e}   {e^{T}Se} = 
\frac{2e^{T}S_{12}e}   {e^{T}Se}
\]
is minimized over all symmetric permutations of $S$. The quantity $\sigma(S,n_{1})$ is called the \emph{uncoupling measure} of $S$ with respect to parameter $n_{1}$. In other words $\sigma(S,n_{1})$ is the ratio of the sum of the elements in the off-diagonal blocks to the sum of all the matrix entries.
\end{definition}

Before moving on, two points should be made clear. First, there is no arbitrary uncoupling measure value below which a matrix is deemed to be nearly uncoupled. Rather, $\sigma(S,n_{1})$ is a relative value whose meaning is dependent on the uncoupling measures of $S$ using other choices of $n_{1}$ or on comparisons with other similarity matrices a researcher has experience with. Second, exact calculation of the uncoupling measure for all but very small problems is not feasible, but its theoretical value is  important since it allows us to compare matrices $S$ and $P$ as the the following theorem shows.

\begin{theorem} \label{th:sigma} 
If $S$ is the $n \times n$ consensus matrix created from $r$ clustering results, then for the doubly stochastic matrix $P=DSD$, $\sigma(P,n_{1}) \le \frac{\Sigma}{nr}\sigma(S,n_{1}) $, where $\Sigma = e^{T}Se$.
\end{theorem}

\begin{proof}
By the way we constructed $S$, $s_{ii}=r$ for $i=1,2,\dots,n$. Since $p_{ii}=d_{i}d_{i}s_{ii}$ and $p_{ii} \le 1$, it follows that $d_{i}^{2}r$ implies  $d_{i} \le \frac{1}{\sqrt{r}}$.

If we impose the same block structure on $D$ that exists for $S$, that is
\[
D = \left[
\begin{array}{cc}
 D_{1} & 0     \\
 0 &  D_{2}          
\end{array}
\right],
\]
and recall that $P$ is doubly stochastic,
\[
\sigma(P,n_{1}) = \frac{2e^{T}D_{1}S_{12}D_{2}e}   {n}.
\]
Since each element of $D_{1}$ and $D_{2}$ is less than $\frac{1}{\sqrt{r}}$,
\[
\sigma(P,n_{1}) \le \frac{\left(\frac{1}{\sqrt{r}}\right)^{2}(2e^{T}S_{12}e)}{n} = \frac {\Sigma} {nr}  \sigma(S,n_{1}), 
\]
\noindent and the bound is established.
\end{proof}

\subsection{The spectrum of $P$} \label{sec:spec}

Consider the following facts about the eigenvalues of $P$.

\begin{enumerate}
\item Since $P$ is stochastic, all of its eigenvalues lie on or inside the unit circle of the complex plane.
\item Since $P$ is real-symmetric, all of its eigenvalues are real. Combined with the last fact, this means all eigenvalues of $P$ reside in the interval $[-1,1]$.
\item The largest eigenvalue of $P$ is one, and since $P$ is irreducible, that eigenvalue is simple (i.e. it appears only once).
\item $\lambda_{i}(P) \ne -1$ for all $i$ because $P$ is a primitive matrix. $P$ is primitive because it is irreducible and has at least one positive diagonal element (\cite{meyer2000maa}, p. 678). 
\end{enumerate} 

Unlike Markov chain researchers who desire a small second eigenvalue since it leads to faster convergence when calculating the chain's stationary distribution vector, we want a second eigenvalue near one. Slow convergence is a good thing for us since it allows time to examine the elements of $x_{t}$ as it passes through short-term stabilization and middle-run evolution. Also, $\lambda_{2}(P) \approx 1$ \emph{may} indicate that the matrix is nearly uncoupled \cite{stewart1994ins}.

We will now show that $\lambda_{2}(P) \approx 1$ along with other properties of $P$ \emph{guarantees} that $P$ is nearly uncoupled. First, observe the following lemma whose proof is self-evident.

\begin{lemma} \label{lem:sym}
Let $\{P_{k}\}$ be a sequence of matrices with limit $P_{0}$. Then,
\begin{enumerate}
\item If each matrix in $\{P_{k}\}$ is symmetric, $P_{0}$ is symmetric, and
\item If each matrix in $\{P_{k}\}$ is stochastic, $P_{0}$ is stochastic.
\end{enumerate}
\end{lemma}




\begin{theorem} \label{thm:hgeasy} 
For a fixed integer $n > 0$, consider the $n \times n$ irreducible, symmetric, doubly stochastic matrix $P$. Given $\epsilon > 0$, there exists a $\delta > 0$ such that if $\sigma(P,n_{1}) < \delta$, then $|\lambda_{2}(P)-1| < \epsilon$. In other words, if $P$ is sufficiently close to being uncoupled, then $\lambda_{2}(P) \approx 1$.
\end{theorem}

\begin{proof}
Two proofs will be presented. The first relies on a continuity argument, while the second gives an explicit bound on $|\lambda_{2}(P)-1|$.

Proof (1): Let $\epsilon > 0$. Consider a sequence of irreducible, symmetric, doubly stochastic matrices 
\[
P_{k}= \left[
\begin{array}{cc}
 P_{11}^{(k)} & P_{12}^{(k)}     \\
 P_{21}^{(k)} & P_{22}^{(k)}  
\end{array}
\right]
\]
defined so that $\displaystyle \lim_{k \rightarrow \infty} \sigma(P_{k},n_{1}) = 0$. The Bolzano-Weierstrass theorem (\cite{bartle1964era}, p. 155) guarantees that this bounded sequence has a convergent subsequence $P_{k_{1}}, P_{k_{2}}, \dots$ which converges to a stochastic matrix $T$ whose structure is
\[
T=
\left[
\begin{array}{cc}
 T_{11} & 0     \\
 0 & T_{22}  
\end{array}
\right], \quad T_{11}\ne0, T_{22}\ne0, 
\]
where each $T_{ii}$ is stochastic. By the continuity of eigenvalues, there exists a positive integer $M$ such that for $k_{i} > M$,
\[
|\lambda_{2}(P_{k_{i}}) - \lambda_{2}(T)| < \epsilon \quad \Rightarrow \quad |\lambda_{2}(P_{k_{i}}) - 1| < \epsilon,
\]
and the theorem is proven.

Proof (2): Suppose that the rows and respective columns have been permuted so that 
\[
P= \left[
\begin{array}{cc}
 P_{11} & P_{12}     \\
 P_{21} & P_{22}  
\end{array}
\right],
\]
where $P$ is nearly uncoupled, and define $C$ to be the $n \times n$ block diagonal matrix with the stochastic complements of $P_{11}$ and $P_{22}$ on the diagonals, that is 
\[
C= \left[
\begin{array}{cc}
 C_{11} & 0     \\
 0 & C_{22}  
\end{array}
\right].
\]
If $E$ is defined to make the equation $C=P+E$ true, then a consequence of the Courant-Fisher Theorem can be used (\cite{meyer2000maa}, pp.~550-552) to show that for any matrix norm\footnote{If the 2-norm is used the bound is $|1-\lambda_{2}(P)| \le 2\sqrt{n}\sigma(P,n_{1})$. We thank Ilse Ipsen for this observation.}
\[
\lambda_{2}(P) - ||E|| \le 1 \le \lambda_{2}(P) + ||E|| \rightarrow  |1-\lambda_{2}(P)| \le  ||E||. 
\] 
\end{proof}
  
\begin{theorem} \label{thm:hghard} 
For a fixed integer $n > 0$, consider the $n \times n$ irreducible, symmetric, doubly stochastic matrix $P$. Given $\epsilon > 0$, there exists a $\delta > 0$ such that if $|\lambda_{2}(P)-1| < \delta$, then $\sigma(P,n_{1}) < \epsilon$ for some positive integer $n_{1}<n$. In other words, if $\lambda_{2}(P)$ is sufficiently close to 1, then $P$ is nearly uncoupled.
\end{theorem}

\begin{proof}
The argument is by contradiction and similar to one used in \cite{hartfiel1998ssm}. Suppose there is an $\epsilon > 0$ such that for any $\delta > 0$ there is an $n \times n$ irreducible, symmetric, doubly stochastic matrix $P$ with $|\lambda_{2}(P)-1| < \delta$ and $\sigma(P,n_{1}) > \epsilon$ for all for positive integers $n_{1}<n$. For $\delta = \frac{1}{k}$ let $P_{k}$ be such a matrix. There must be a subsequence $P_{i_{1}}, P_{i_{2}}, \dots$ which converges, say to $P_{0}$. Then $P_{0}$ must have  $\lambda_{2}(P_{0}) = 1$ and thus $\sigma(P_{0},n_{1}) = 0$. Yet, $\sigma(P_{0},n_{1}) = \lim_{k \rightarrow \infty} \sigma(P_{k},n_{1}) \ge \epsilon$, a contradiction. 
\end{proof}

Although we previously defined an uncoupling measure for a general matrix in Section \ref{sec:dsd}, for doubly stochastic matrices this theorem allows us to use $\lambda_{2}$ as an uncoupling indicator with a value near one signifying almost complete uncoupling.

There may be additional eigenvalues of $P$ that are close to one. This group of eigenvalues is called the \emph{Perron cluster} \cite{deuflhard2000iai, deuflhard2005rbc}, and in the case where all eigenvalues are real the Perron cluster can be defined as follows.

\begin{definition}
Let $P$ be an $n \times n$ symmetric, stochastic matrix with eigenvalues, including multiplicities, of $1=\lambda_{1} \ge \lambda_{2} \ge \lambda_{3} \ge \dots \ge \lambda_{n}$. If the largest difference between consecutive eigenvalues occurs between $\lambda_{k}$ and $\lambda_{k+1}$, the set $\{1, \dots \lambda_{k}\}$ is called the \emph{Perron cluster} of $P$. If two or more pairs of eigenvalues each have differences equal to the largest gap, use the smallest value of $k$ to choose $\lambda_{k}$.  The larger the gap, the more well-defined the cluster.
\end{definition}

Some researchers use the number of eigenvalues in the Perron cluster as the number of clusters they search for \cite{deuflhard2000iai, fritzsche2008sai}. This inference is a natural extension of Theorems \ref{thm:hgeasy} and \ref{thm:hghard}, that is if $P$ had $k$ eigenvalues sufficiently close to 1, then $P$ is nearly uncoupled with $k$ dominant diagonal blocks emerging after an appropriate permutation $QPQ^{T}$. This is also the approach we will take with the stochastic clustering algorithm. Unlike with the vast majority of clustering methods, the user will not have to tell the algorithm the number of clusters in the data set unless they explicitly want to override the algorithm's choice. Instead, the stochastic clustering algorithm will set $k$ equal to the size of the Perron cluster.

\section{Putting the concept into practice} \label{sec:cip}

Now that the theoretical underpinnings are in place, it is time to formally describe the stochastic clustering algorithm. 

The algorithm takes as input the consensus similarity matrix $S$ which the user has created from whatever combination of clustering methods and/or parameter settings they choose. $S$ is then converted into the doubly stochastic matrix $P$ using the Sinkhorn-Knopp algorithm.  All eigenvalues are computed, and the Perron cluster of $P$ is identified. Eigenvalues of symmetric matrices can be efficiently computed  \cite{parlett1998sep}, but if finding all eigenvalues is too costly, the user, with knowledge of the underlying data set, can direct the program to find only the $\hat{k}$ largest eigenvalues \textcolor{blue}{($\hat{k}>k$)}. The size, $k$, of the Perron cluster of these $\hat{k}$ eigenvalues is then used by the stochastic clustering algorithm to separate the data into $k$ clusters.

Starting with a randomly generated $x_{0}^{T}$, $x_{t}^{T}=x_{t-1}^{T}P$ is evaluated for $t=1,2,\dots$. After each calculation, the entries of $x_{t}^{T}$ are sorted, the $k-1$ largest gaps in the sorted list identified and used to divide the entries into $k$ clusters.  When the $k$ clusters have been identical for $n$ iterations, where $n$ is a user-chosen parameter, the program stops and the clusters returned as output. Figure \ref{fig:sca} summarizes the algorithm.

\begin{figure}
\centering
\headerframe {5pt}{.5pt}{
\begin{minipage}[t]{5in}
\textbf{Stochastic Clustering Algorithm (SCA)}
\begin{enumerate}
\item Create the consensus similarity matrix $S$ using a clustering ensemble of user's choice.
\item Use matrix balancing to convert $S$ into a doubly stochastic symmetric matrix $P$.
\item Calculate the eigenvalues of $P$. The number of clusters, $k$, is the number of eigenvalues in the Perron cluster.
\item Create a random $x_{0}^{T}$.
\item Track the evolution $x_{t}^{T}=x_{t-1}^{T}P$. After each multiplication, sort the the elements of $x_{t}^{T}$ and then separate the elements into $k$ clusters by dividing the sorted list at the $k-1$ largest gaps. Alternatively, the elements of $x_{t}$ can be clustered using $k$-means or any other widely available clustering method. When this clustering has remained the same for a user-defined number of iterations, the final clusters have been determined. 
\end{enumerate}
\end{minipage}
}
\caption{The Stochastic Clustering Algorithm}
\label{fig:sca}
\end{figure}

\subsection{A Small Example}
Consider the following small data matrix which includes the career totals in nine statistics for six famous baseball players (the row labels stand for Games, Runs, Hits, Doubles, Triples, Home Runs, Runs Batted In, Stolen Bases, and Bases on Balls).
\[
A=\bordermatrix{  & \mbox{Rose} & \mbox{Cobb} & \mbox{Fisk} & \mbox{Ott} & \mbox{Ruth} & \mbox{Mays} \cr
\mbox{G} &  \hfill 3562 & \hfill 3034 & \hfill 2499 & \hfill 2730 & \hfill 2503 & \hfill 2992 \cr
\mbox{R} &  \hfill 2165 & \hfill 2246 & \hfill 1276 & \hfill 1859 & \hfill 2174 & \hfill 2062 \cr
\mbox{H} & \hfill 4256 & \hfill 4189 & \hfill 2356 & \hfill 2876 & \hfill 2873 & \hfill 3283 \cr
\mbox{2B} &  \hfill 746 & \hfill 724 & \hfill 421 & \hfill 488 & \hfill 506 & \hfill 523 \cr
\mbox{3B} &  \hfill 135 & \hfill 295 & \hfill 47 & \hfill 72 & \hfill 136 & \hfill 140 \cr
\mbox{HR} &  \hfill 160 & \hfill 117 & \hfill 376 & \hfill 511 & \hfill 714 & \hfill 660 \cr
\mbox{RBI} & \hfill 1314 & \hfill 1938 & \hfill 1330 & \hfill 1860 &  \hfill 2213 & \hfill 1903 \cr
\mbox{SB} & \hfill 198 & \hfill 897 &\hfill 128 & \hfill 89  & \hfill 123 & \hfill 338 \cr
\mbox{BB} & \hfill 1566 & \hfill 1249 & \hfill 849 & \hfill 1708 &  \hfill 2062 & \hfill 1464 } .
\]

Those familiar with baseball history would mentally cluster these players into singles hitters (Rose and Cobb), power hitters (Mays, Ott, and Ruth) and, a great catcher who doesn't have enough home runs and runs batted in to fit with the power hitters nor the long career and large number of hits to fit with the singles hitters (Fisk). 

The consensus similarity matrix was built using  the multiplicative update version of the nonnegative matrix factorization algorithm \cite{leeseung}. Since it isn't clear whether two or three clusters would be most appropriate, $S$ was created by running this algorithm 50 times with $k=2$ and 50 times with $k=3$. The resulting similarity matrix is
\[
S=\bordermatrix{ & \mbox{Rose} & \mbox{Cobb} & \mbox{Fisk} & \mbox{Ott} & \mbox{Ruth} & \mbox{Mays} \cr
\mbox{Rose} &  \hfill 100 & \hfill 67 & \hfill 73 & \hfill 2 & \hfill 0 & \hfill 2 \cr
\mbox{Cobb} &  \hfill 67 & \hfill 100 & \hfill 50 & \hfill 1 & \hfill 2 & \hfill 7 \cr
\mbox{Fisk} & \hfill 73 & \hfill 50 & \hfill 100 & \hfill 15 & \hfill 9 & \hfill 24 \cr
\mbox{Ott} &  \hfill 2 & \hfill 1 & \hfill 15 & \hfill 100 & \hfill 92 & \hfill 82 \cr
\mbox{Ruth} &  \hfill 0 & \hfill 2 & \hfill 9 & \hfill 92 & \hfill 100 & \hfill 77 \cr
\mbox{Mays} &  \hfill 2 & \hfill 7 & \hfill 24 & \hfill 82 & \hfill 77 & \hfill 100  } .
\]

With a small example like this, especially one where the players that will cluster together have been purposely placed in adjacent columns, it would be simple enough to cluster the players through a quick scan of $S$. However, following the algorithm to the letter we apply the Sinkhorn-Knopp algorithm. The resulting doubly stochastic matrix (rounded to four places) is 

\[
P=\bordermatrix{ & \mbox{Rose} & \mbox{Cobb} & \mbox{Fisk} & \mbox{Ott} & \mbox{Ruth} & \mbox{Mays} \cr
\mbox{Rose} & 0.4131 &  0.2935  &  0.2786  &   0.0075   &      0        &  0.0075\cr
\mbox{Cobb} &0.2935   & 0.4644  &   0.2023 &   0.0040   &   0.0082 &  0.0277\cr
\mbox{Fisk}   &0.2786   &   0.2023 &  0.3525  &   0.0517   &   0.0323 &  0.0826\cr
\mbox{Ott}      &0.0075   &  0.0040  &  0.0517  &  0.3374    &  0.3233  &  0.2761\cr
\mbox{Ruth}  &0             &  0.01082 &  0.0323  &   0.3233  &   0.3660  & 0.2701\cr
\mbox{Mays}  &0.0075   &  0.0277  &  0.0826  &   0.2761  &   0.2701  &  0.3361 } .
\]

The eigenvalues of $P$ are 1.0000,  0.8670,  0.2078, 0.1095, 0.0598,  and 0.0254 suggesting that there are two clusters in this data.

Table \ref{baseballtable} shows the results from a sample run of the clustering method. The initial probability vector  $x_{0}^{T}$ was chosen randomly, and the table shows the value of $x_{t}^{T}$ and the corresponding clusters for the next seven steps of the algorithm. Since $k=2$, the clusters are determined by ordering the entries of $x_{t}^{T}$, finding the largest gap in this list, and clustering the elements on either side of this gap. For example, at the $t=6$ step shown in the table, the largest gap in the sorted list is between $0.1609$ and $0.1715$. This leads to the numerical clustering of $\{ 0.1597, 0.1600, 0.1609\}$ and $\{0.1715, 0.1739, 0.1741\}$ which translates to the clustering \{Rose, Cobb, Fisk\} and \{Ott, Ruth, Mays\}.

\begin{table}[h]
\caption{Following the Stochastic Clustering Algorithm for the Small Example}
\begin{footnotesize}
\begin{center}
\begin{tabular}{c | c | c}
$t$ & $x_{t}^{T}$ & Clusters\\ \hline
0 & 
$\left( 
\begin{array}{rrrrrr}
 0.2334  &  0.2595   & 0.0364  &  0.2617  &  0.1812  &  0.0279
\end{array}
\right)$ 
&  \{Rose, Cobb, Ott, Ruth\}  \\
& & \{Fisk, Mays\} \\
\hline
1 &
$\left( 
\begin{array}{rrrrrrrrr}
 0.1848 &   0.1997  &  0.1520  &  0.1592 &   0.1618  &  0.1425
\end{array}
\right)$
& \{Rose, Cobb\} \\
& & \{Fisk, Ott, Ruth, Mays\} \\
\hline
2 &
$\left( 
\begin{array}{rrrrrrrrr}
0.1795 &   0.1836 &   0.1707 &   0.1554&    0.1557  &  0.1550
\end{array}
\right)$
& \{Rose, Cobb, Fisk\} \\
& & \{Ott, Ruth, Mays\} \\
\hline
3 &
$\left( 
\begin{array}{rrrrrrrrr}
0.1779  &  0.1787 &   0.1732 &  0.1565  &  0.1561 &   0.1576
\end{array}
\right)$
& \{Rose, Cobb, Fisk\} \\
& & \{Ott, Ruth, Mays\} \\
\hline
4 &
$\left( 
\begin{array}{rrrrrrrrr}
0.1765  &  0.1765 &   0.1729 &   0.1578  &  0.1574 &   0.1589
\end{array}
\right)$
& \{Rose, Cobb, Fisk\} \\
& & \{Ott, Ruth, Mays\} \\
\hline
5 &
$\left( 
\begin{array}{rrrrrrrrr}
0.1752 &   0.1751 &   0.1722  &  0.1590 &   0.1586 &   0.1600
\end{array}
\right)$
& \{Rose, Cobb, Fisk\} \\
& & \{Ott, Ruth, Mays\} \\
\hline
6 &
$\left( 
\begin{array}{rrrrrrrrr}
0.1741  &  0.1739  &  0.1715  &  0.1600  &  0.1597 &   0.1609
\end{array}
\right)$
& \{Rose, Cobb, Fisk\} \\
& & \{Ott, Ruth, Mays\} \\
\hline
7 &
$\left( 
\begin{array}{rrrrrrrrr}
0.1731 &   0.1729  &  0.1709 &   0.1609  &  0.1606 &   0.1616
\end{array}
\right)$
& \{Rose, Cobb, Fisk\} \\
& & \{Ott, Ruth, Mays\} \\
\hline
\end{tabular}
\end{center}
\end{footnotesize}
\label{baseballtable}
\end{table}

From $t=2$ on the clusters remain the same. The SCA defines the stopping condition as a user-defined number of consecutive identical clusterings. If that number is six, then the final clustering of \{Rose, Cobb, Fisk\} and \{Ott, Ruth, Mays\} is determined when $t=7$.  For the reader wondering if the clustering changes at some later point, the algorithm was run through $t=1000$ and the same clustering was found at each step.

\section{Implementation}
As is to be expected with a new algorithm, actual implementation of ideas that looked fine on paper can still be problematic. Even before implementation, there may be concerns about perceived weak links in the algorithm. In this section we will address some of these concerns. Since this section and the results section involve many of the same issues, it will be hard to talk about them without highlighting some of the results to come. Hopefully, no great surprises are spoiled, and the turning of pages back and forth is kept to a minimum.

\subsection{Impact of initial probability vectors}  \label{sec:ipvds}

The fact that the stochastic clustering algorithm depends on a random initial probability vector (IPV) raises the question of whether all random probability vectors will lead to the same clustering. Since $P$ is irreducible, we are guaranteed that the matrix has a unique stationary distribution vector that is independent of the IPV. But, for clustering purposes, that is not the issue. Instead we would like to have confidence that for a certain IPV, $x_{t}^{T}$ will remain in short-term stabilization and middle-run evolution long enough for us to identify the clusters. Secondly, as we will see soon in Section \ref{sec:res}, different IPVs can lead to different cluster results.

We will consider the IPV question in two parts. First we address the rare occurrence of an IPV that does not lead to a clustering at all, and then we address the fact that different IPVs can lead to different clusterings.

\subsection{IPVs leading to no solution}
Clearly not every initial probability vector will help us in data clustering. Suppose, for example, that 
$$
x_{0}^{T} = \left( \frac{1}{n} \quad  \frac{1}{n} \quad \frac{1}{n} \quad \dots  \quad \frac{1}{n}  \right)_{1 \times n}.
$$
Since $P_{n \times n}$ is doubly stochastic, $x_{0}^{T}$ is its stationary distribution vector. With such a choice for the IPV, $x_{t}^{T}$ never changes and we have no ability to group the probabilities in $x_{t}^{T}$ in order to cluster the original data.

It is simple enough to make sure that $x_{0}^{T}$ is not the uniform distribution vector, but it is equally important that there are enough iterations for the algorithm to recognize either short-term stabilization or middle-run evolution before $x_{t}^{T}$ reaches the uniform vector. Since each new $x_{t}^{T}$ is the result of the continuous operation of matrix multiplication, $x_{t}^{T}$ being close to the uniform distribution vector, ensures that $x_{t+1}^{T}$ can not be significantly further away for it. Therefore, even though the algorithm generates $x_{0}^{T}$ randomly, the cautious user may want to set a tolerance $\epsilon$ and if
\[
|| x_{0}^{T} - (1/n \;\; 1/n \;\; \dots \;\; 1/n) || < \epsilon,
\]
generate another $x_{0}^{T}$. It should be noted that in the preparation of this paper the stochastic clustering algorithm was run hundreds, if not thousands, of times and never was a failure due to an IPV being too close to the uniform distribution.

\subsection{IPVs leading to different solutions}

The fact that cluster analysis is an exploratory tool means that getting different solutions depending on the initial probability vector is not the end of the road, but rather an opportunity to examine these solutions in the hope of gaining additional insight into the data set's structure.

That said, it would still be instructive to know as much as possible about the characteristics shared by IPVs that lead to the same solution, how many different solutions are possible, and how often each of them is likely to appear. Probabilistic analysis of random starting vectors has been done in the context of iterative methods for finding eigenvalues and eigenvectors \cite{dixon1983eee, jessup1992iai}, and is a natural area for further research on the stochastic clustering method.

\subsection{Using a single measure} \label{sec:singlemeasure}
The workload in consensus clustering is concentrated at the beginning of the process when the large number of clustering results are computed. Even if a user has access to a multiprocessor environment where this work can be shared, it would be advantageous to find a single similarity measure which is compatible with the stochastic clustering algorithm.

\begin{figure}[ht]
\centering
\includegraphics[scale=.4,trim=0 0 0 0]{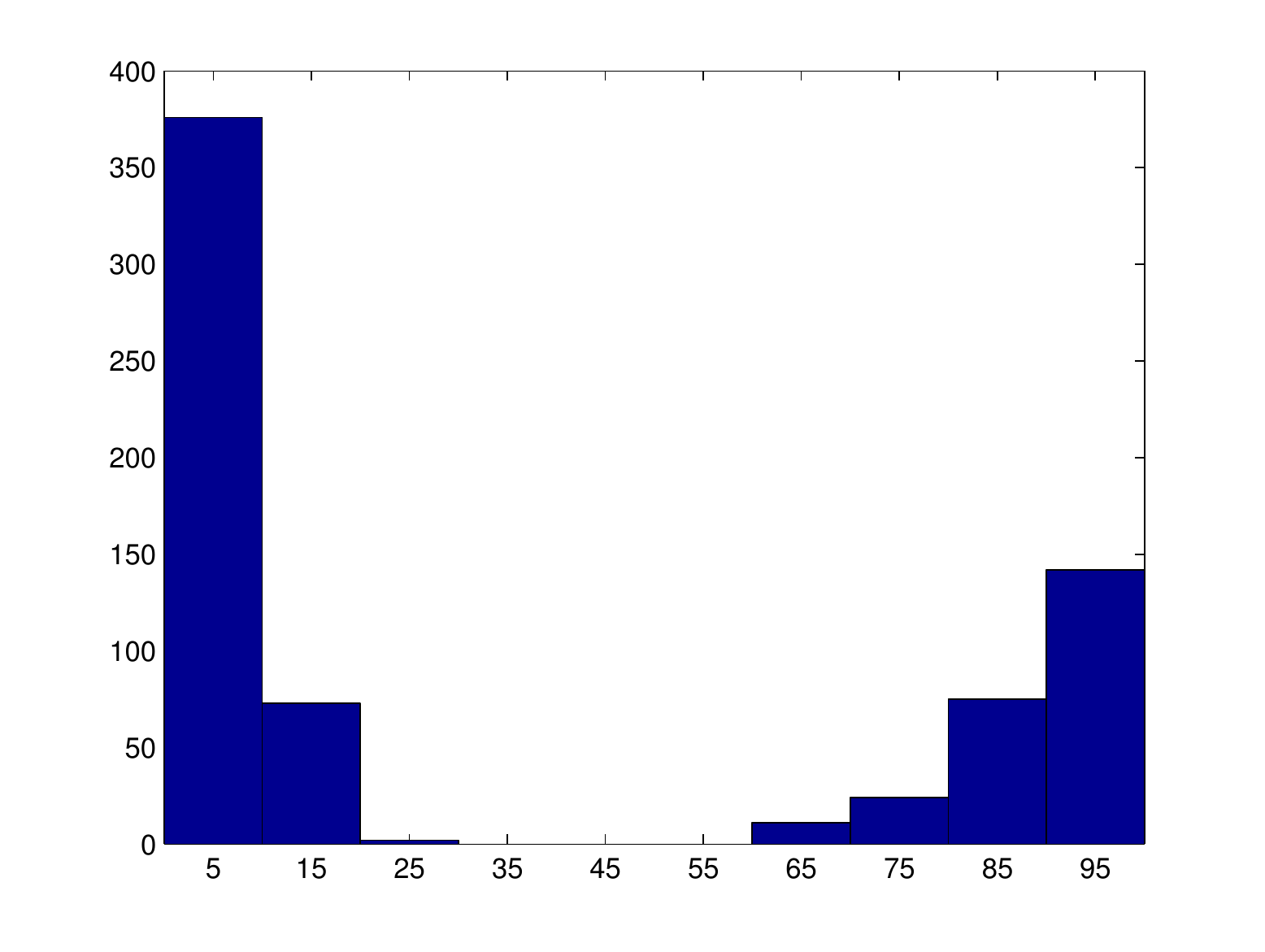}
\caption{This is the histogram of the 703 similarity values used to build a consensus matrix for the 38-element leukemia DNA microarray data set that will be introduced in Section \ref{sec:dna}. The horizontal axis measures the number of times out of 100 that two elements clustered together. The histogram shows that pairs of data points clustered together either a small or large number of times.}
\label{fig:scm}
\end{figure}

\begin{figure}[ht]
\includegraphics[scale=.4,trim=40 0 0 225]{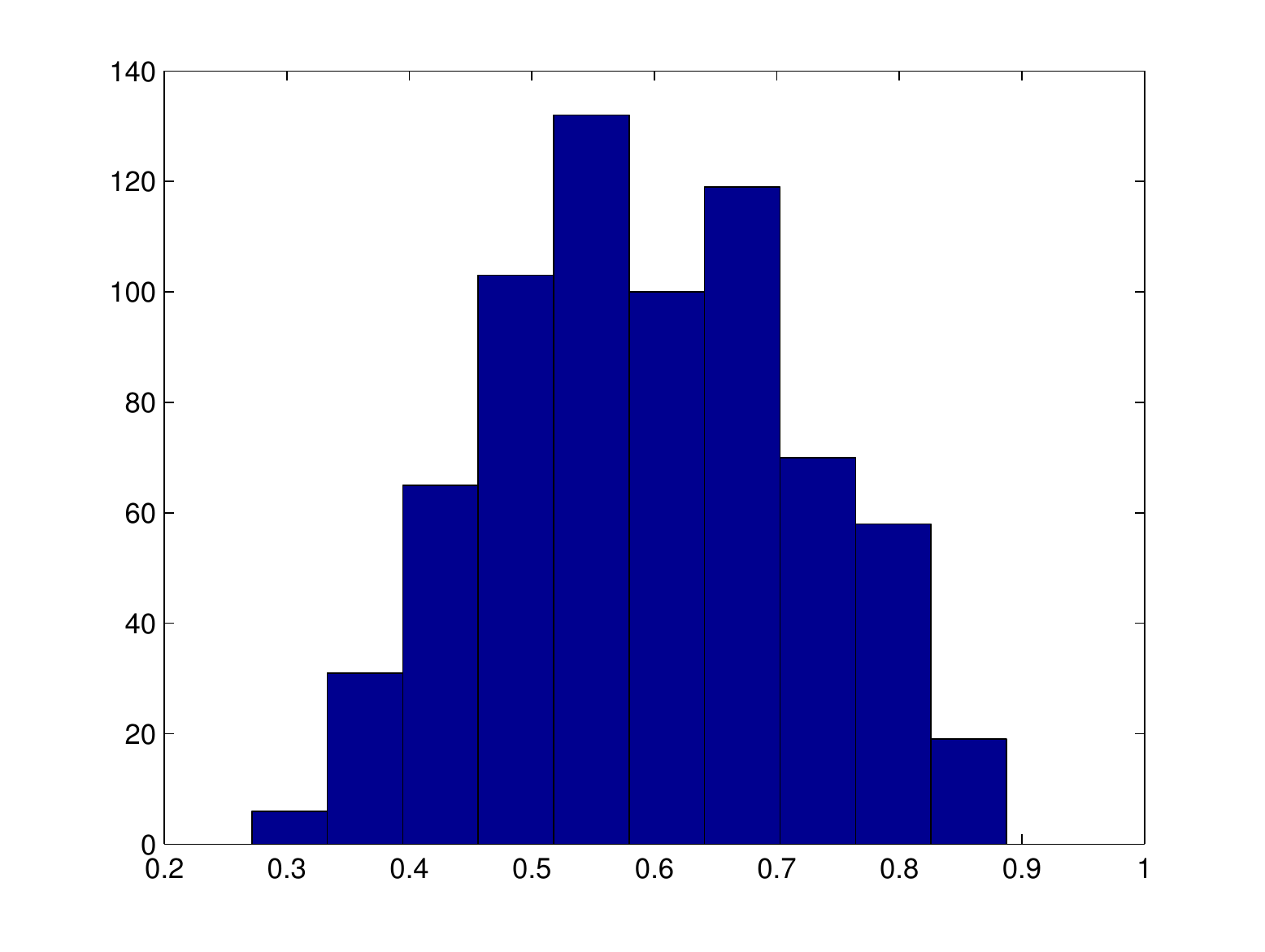}
\includegraphics[scale=.4,trim=40 225 0 225]{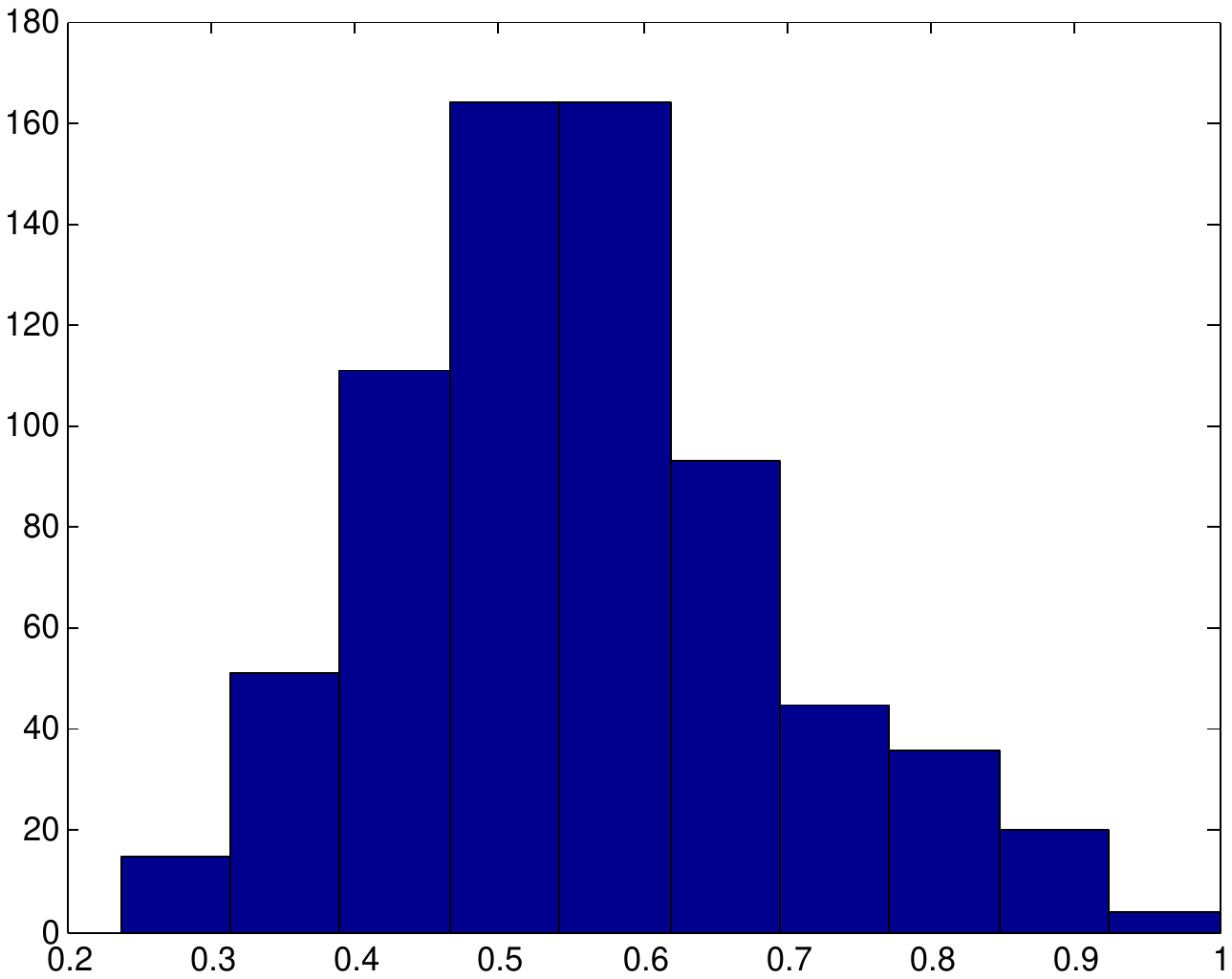}
\caption{The histogram on the left shows the distribution of cosine similarity measures between the same elements used for Figure \ref{fig:scm}, while the histogram on the right does the same for Euclidean norm values scaled to the interval $[0,1]$. Contrast these distributions with the one shown in Figure \ref{fig:scm}.}
\label{fig:simvals}
\end{figure}

Since the SCA is inspired by Simon-Ando theory, the underlying matrix must be nearly uncoupled. For a given data set, the problem with most traditional similarity (or dissimilarity) measures is that their values tend to the middle of their range. To illustrate, consider two common similarity measures: Euclidean distance and the cosine measure
$$
c({x}_{1},{x}_{2}) = \frac{{x}_{1}^{T} {x}_{2}} {||{x}_{1}||_{2} \, ||{x}_{2}||_{2}}.
$$
The former has the advantage of being familiar to almost everyone, while the latter has been found to be particularly useful in text-mining \cite{berry}. However, as Figures \ref{fig:scm} and \ref{fig:simvals} show for the leukemia DNA microarray data set that will be introduced in Section \ref{sec:dna}, the distribution of values returned by these two common measures is not the kind of distribution needed to form a nearly uncoupled matrix. 

In the case of the cosine measure whose range is $[0,1]$, there have been attempts to ``massage'' distributions so that they contain more values near the extremes. Such methods often involve changing small values to zero and then performing some arithmetic operation that gives the remaining data a larger variance (for example, squaring each value) \cite{vandongen2000gcf}. These methods, however, are far from subtle and in experiments for use with the SCA, the matrix $P$ went from dense to too sparse for clustering in one iteration of attempting to adjust its values.

Working with the Euclidean norm brings with it the additional requirement large distances need to be mapped to small similarity values while small distances are mapped to large similarity values. A typical function used in making this translation is a Gaussian of the form
\[
f({x}_{1},{x}_{2}) = e^{\frac{-||x_{1}-x_{2}||^{2}}{2\sigma^{2}}},
\]
where $\sigma$ is a parameter that typically has to be adjusted for each similarity matrix \cite{mouysset2008hpc}. This is certainly an area for future study in implementing the SCA, but so far a reliable way to build a matrix of Gaussians with the distribution required by the SCA has not been found.

It should be noted that \emph{power iteration clustering} introduced by Lin and Cohen has succeeded in using a single measure to cluster data using an algorithm similar in philosophy to the SCA. This method uses a row-stochastic Laplacian-like matrix derived from a similarity matrix constructed using the cosine similarity measure \cite{lin2010pic, lin2010vfm, meila2001rwv}. Like the SCA, clusters are determined by examining intermediate iterates of the power method. It is interesting to note despite mentioning a Gaussian approach to the Euclidean norm in \cite{lin2010pic}, all results in the paper were obtained using either a $0-1$ or cosine measure.

A single measure that has been used with some success involves the idea of nearest neighbors, those data points closest to a given data point using a specific distance measure. For each element $g$ in the data set, let the set $\mathcal{N}_{g}$ consist of the $\kappa$ nearest neighbors of $g$, where the user chooses both the positive integer $\kappa$ and the distance measure used. The $s_{ij}$ element of the consensus matrix is equal to the number of elements in $\mathcal{N}_{i} \cup \mathcal{N}_{j}$ \cite{abbey2011pc}.

Work with consensus matrices built in this fashion is still in its initial stages. It has become obvious that the choice of $\kappa$ and the distance measure greatly affect the results as can be seen in Table \ref{tbl:knn}.

\begin{table}[ht]
\caption{Building a consensus matrix based on the number of shared nearest neighbors can work well or poorly depending on the value of $\kappa$, the number of nearest neighbors calculated for each data point. The results in this table are from clustering the rather simple, four-cluster Ruspini data set \cite{ruspini}. When $\kappa=15$ the stochastic clustering algorithm detects five clusters. This fifth cluster only has one member, while the rest of the solution is correct.}
\begin{center}
\begin{tabular}{|c|c|c|} \hline
$\kappa$ & Clusters & Errors \\ \hline
15 & 5 & 1 \\
20 & 4 & 0 \\
25 & 4 & 18 \\ \hline
\end{tabular}
\end{center}
\label{tbl:knn}
\end{table}

 
\section{Results} \label{sec:res}
In building test cases for our proposed algorithm, one complication is determining the ensemble used to build the initial similarity matrix $S$. In the results that follow the ensembles will typically consist of multiple runs of multiplicative update version of NMF \cite{leeseung} or $k$-means or a combination of both.\footnote{For an example of how the factors found by NMF are used to cluster data see \cite{chartier2011nap}.} In each case, the value or values of $k$ used when calling these algorithms will be noted, though as explained above the new stochastic clustering algorithm will use the number of eigenvalues in the Perron cluster of $P$ to determine $k$.

\subsection{Iris data set}
The Fisher iris data set \cite{fisher} consists of four measurements (petal length, petal width, sepal length, and sepal width) for 150 iris flowers, fifty each from three iris species (\emph{Iris setosa}, \emph{Iris virginica}, and \emph{Iris versicolor}). It is well-documented that the \emph{setosa} cluster is easily separable from the other two, but separating the \emph{virginica} and \emph{versicolor} species is more difficult \cite{fredjain2002}.

When building $S$ using NMF the choice of $k$ is limited to two or three since NMF requires $k$ to be less than both the dimension of the data and the number of samples. Running the multiplicative update version of NMF 100 times with $k=2$ never results in a perfect splitting of \emph{setosa} from the other two species, though there are three or fewer clustering errors 67 times. However, there are six instances of more than 15 errors including a worst case of 26. Despite these problems, the SCA, using a consensus similarity matrix built from these rather poor results gets the clustering correct for all but three irises. Although NMF does quite poorly in trying to separate the irises into three clusters, the $S$ derived from these results leads to a perfect two-cluster separation of \emph{setosa} irises from \emph{virginica} and \emph{versicolor} ones. 

On the whole, individual clustering results on the iris data set using $k$-means clustering with $k=2$ or $k=3$ are better than those returned by NMF. However, building $S$ using the results from $k$-means clustering, we get very similar results to what we saw with NMF.

If we decide to build $S$ using $k=4$ just to see if it will give us any insight into the data set, SCA recognizes that there are three clusters in the data set, but 16 flowers are misclustered. Though that result may not seem encouraging, notice that this an improvement over the the range of errors (21 - 38) when using $k$-means with $k=3$. 


Finally, the consensus matrices found using NMF and $k$-means were summed to see if a more robust clustering than the one found by SCA using $S$ from just one of these methods could be found. Notice that this approach proved fruitful as there is at most one error regardless of the value of $k$ used.

The results from using all of these different consensus matrices are summarized in Table \ref{irisresults1}.

\begin{table}[h!]
\centering
\caption{\noindent Clustering the iris data set, $S$ created using NMF (first two lines), $k$-means (next three lines), and a combination of the two (last two lines).} 
\begin{tabular}{c | c | c | c}
Method and $k$  & range of \#  of errors &      $k$ found & \# of errors    \\ 
  used  to create $S$    & in single clusterings & by SCA & in SCA result \\ \hline \hline
NMF (2) & 1--26 & 2 & 3  \\
NMF (3) & 19--72 & 2 & 0    \\ \hline 
$k$-means (2) & 3 & 2 & 3  \\
$k$-means (3) & 21--38 & 2 & 0   \\
$k$-means (4) & n/a & 3 & 16  \\ \hline 
Combined (2) & 1--26 & 2  & 1   \\ 
Combined (3) & 19--72 & 2  & 0   \\ 
\end{tabular}
\label{irisresults1}
\end{table}  

\subsection{Leukemia DNA microarray data set} \label{sec:dna}
In 1999 a paper was published analyzing a DNA microarray data set containing the gene expression values for 6817 genes from 38 bone marrow samples \cite{leukemia}. Five years later, the same 38 samples were examined, though this time only 5000 genes were used \cite{metagene}. The samples came from leukemia patients who had all been diagnosed with either acute lymphoblastic leukemia (ALL) or acute myeloid leukemia (AML). Additionally, the ALL patients had either the B-cell or T-cell subtype of the disease (ALL-B or ALL-T). This data set is well known in the academic community (Google Scholar reports that the 1999 paper has been cited over 6000 times) and is an excellent test for new clustering algorithms since it can be divided into either two (ALL/AML) or three (ALL-B/ALL-T/AML) clusters. The actual clustering for the leukemia data set is known (see Table \ref{clusters}), though the 2004 paper noted that the data ``contains two ALL samples that are consistently misclassified or classified with low confidence with most methods.  There are a number of possible explanations for this, including incorrect diagnosis of the samples  \cite{metagene}.''

\begin{table}[ht]
\caption{The correct clustering of the leukemia DNA microarray data set.}
\begin{center}
\begin{tabular}{|c|c|} \hline
 Diagnosis & Patients  \\ \hline
 ALL-B & 1 -- 19  \\
ALL-T  & 20 -- 27  \\
AML &  28 -- 38   \\ \hline
\end{tabular}
\end{center}
\label{clusters}
\end{table}

Since the 2004 paper was published to demonstrate the effectiveness of nonnegative matrix factorization in clustering this data set, this seems to be an appropriate test for the stochastic clustering algorithm, using NMF with different $k$ values to build the ensemble. The data set was clustered using NMF 100 times each for $k=2$ and $k=3$. Additionally, to explore the data set further, the data were clustered an additional 100 times for $k=4,5$ and 6.

Figure \ref{tbl:dnaerr2} shows the number of errors for each clustering used in building $S_{2}$, the $k=2$ consensus similarity matrix. NMF is clearly quite good at clustering this data set into two clusters, which was the point of \cite{metagene}. Each time the stochastic clustering algorithm is used to cluster the patients based on $S_{2}$, it  makes exactly two errors -- misclustering Patients 6 and 29. 

Similar comparisons were done using $S_{3}$, the $k=3$ consensus similarity matrix, and again the stochastic clustering method could not improve on the already excellent results of NMF. 
NMF made an average of 3.18 errors per clustering compared to 4.76 for the SCA. Even the hope that the SCA would provide a narrower band of errors than NMF is not realized (see Table \ref{tbl:dnaerr3}). Perhaps the lesson is that if the original method does a good job of clustering, then SCA is not likely to improve on it, though it is also not likely to worsen it. 

Since cluster analysis is an exploratory tool, consensus matrices $S_{4}$, $S_{5}$, and $S_{6}$ were constructed to see if either the stochastic clustering algorithm or nonnegative matrix factorization could discover some hidden structure in the data set that would indicate one or more undiscovered clusters. If a group of elements all break away from an existing cluster or clusters, there is reason for further investigation regarding a new cluster. Interestingly, when $k=4$, the results from both NMF and the SCA agree. As Table \ref{tbl:dnaerr4} summarizes, they both have identified a fourth cluster made up of four ALL-B patients and two AML patients.

Neither of the methods give any indication of further clusters. When $k=5$ or $k=6$ both methods begin to build two or three large clusters with the remaining clusters containing only two or three members.

\begin{figure}
\centering
\subfloat[Clustering comparisons when $k=2$][The leukemia DNA microarray data set was clustered 100 times using NMF with $k=2$. The number of errors ranged between one and four. When the SCA was used on the consensus matrix created from those 100 NMF clusterings, it mis-clustered Patients 6 and 29 each time.]
{
\begin{tabular}{ccccc} \hline
\# of Errors & 1 & 2 & 3 & 4  \\ \hline
\# of Instances (NMF) & 30 & 65 & 3 & 2  \\ \hline 
\# of Instances (SCA)  & 0 & 100 & 0 & 0  \\ \hline 
\end{tabular}
\label{tbl:dnaerr2}
} \\ 
\subfloat[Clustering comparisons when $k=3$][Neither the SCA nor NMF shows an advantage over the other when clustering the consensus matrix $S_{3}$. ]
{
\begin{tabular}{ccccccccccc} \hline
\# of Errors & 1 & 2 & 3 & 4 & 5 & 6 & 7 & 8 & 9 & 10+  \\ \hline
\# of Instances (NMF) & 0 & 71 & 3 & 9 & 3 & 3 & 1 & 2 & 0 & 8 \\ \hline 
\# of Instances (SCA) & 0 & 67 & 0 & 0 & 0 & 0 & 0 & 0 & 0 & 33 \\ \hline 
\end{tabular}
\label{tbl:dnaerr3}
} \\ 
\subfloat[A new cluster?][Both NMF and SCA agree that there may be a new cluster. The third column shows the membership of this new cluster and the patients remaining in the other three.]
{
\begin{tabular}{|c|c|c|} \hline
 Diagnosis & Patients & Patients  \\ \hline
 ALL-B & 1 -- 19 & 1, 3, 5, 7 -- 9, 11 -- 14, 16 -- 18 \\ \hline
ALL-T  & 20 -- 27 & 10, 20 -- 27 \\ \hline
AML &  28 -- 38  &  28, 30 -- 35, 37, 38 \\ \hline
New Cluster & & 4, 6, 19, 29, 36 \\ \hline
\end{tabular}
\label{tbl:dnaerr4}
} 
\caption[Summary of results for the leukemia data set]{This is a collection of tables that compare the results of clustering consensus matrices constructed using different $k$-values. The consensus matrices were clustered by both the SCA and NMF. Table \ref{tbl:dnaerr2} compares the results for $k=2$. Table \ref{tbl:dnaerr3} shows very little difference between the two methods when $k=3$. Table \ref{tbl:dnaerr4} shows a possible fourth cluster suggested by both NMF and SCA.}
\label{fig:dnasum}
\end{figure}

Before we move on to the next data set, there is one other interesting result to report. If the stochastic clustering algorithm is run using the sum of $S_{2}$ and $S_{3}$ it identifies two clusters and makes only one clustering mistake, namely Patient 29.\footnote{Throughout the research period for this paper, the Patient 29 sample was misclustered nearly 100 per cent of the time. One of the authors of the 2004 paper verifies that in their work, the Patient 29 sample was also often placed in the wrong cluster \cite{tamayo2011pc}.}

\subsection{Custom clustering}
As we first mentioned in Section \ref{sec:ipvds}, the fact that the stochastic clustering algorithm uses a random initial probability vector means that it can arrive at different solutions, and when clustering the leukemia data set we found this to be so. While this might be viewed as a weakness of the algorithm, it does give the researcher the ability to answer a very specific question by creating a specific initial probability vector. 

In Section \ref{sec:dna}, we noticed that the SCA did not cluster the leukemia data set consensus matrix any better than nonnegative matrix factorization. But what if our primary interest was not in clustering the entire data set, but instead in finding the membership of the cluster of a particular data point. For example, if you are the physician for Patient number 2 you have limited interest in a global view of the leukemia data set. Indeed, rather than knowing which of the three clusters Patient 2 belonged to, it would be of greater use to you to know a small number of other patients that are most like Patient 2 in the hope that that knowledge would help you tailor the best treatment plan possible. 

To create such a custom clustering, we construct an IPV containing all zeros except for a 1 in the place corresponding to our data point of interest. We then ask the stochastic clustering algorithm to find the cluster containing our specific data point. Since we may be interested in a collection much smaller than that cluster, the stochastic clustering algorithm can be modified to ask for a small number data points whose $x_{t}$ entries are closest to our target point.

Here again we find hope in a feature of the SCA that seemed to disappoint us in Section \ref{sec:dna}. In that section, the clustering of consensus matrices built from methods using $k=5$ and $k=6$ seemed to supply new information. In fact, the small clusters found then are indicative of an especially close relationship between the cluster members.

\begin{figure}
\centering
\headerframe{5pt}{.5pt}{
\begin{minipage}[t]{5in}
\textbf{Custom Clustering Algorithm (CCA)}
\begin{enumerate}
\item Create the consensus similarity matrix $S$ and the doubly stochastic symmetric matrix $P$ just as in the stochastic clustering algorithm.
\item Construct $x_{0}^{T}$ to contain all zeros except for a one in the place of the element we are interested in creating a custom cluster for.
\item Pass the algorithm values for the minimum and maximum size cluster you desire and the maximum number of iterations the CCA should take trying to find that cluster.
\item After each $x_{t}^{T}=x_{t}^{T}P$ multiplication, cluster the elements of $x_{t}^{T}$ as in the SCA. If the cluster containing the target element is within the size parameters, output the cluster and end the program.
\end{enumerate}
\end{minipage}
}
\caption{The Custom Clustering Algorithm}
\label{fig:cca}
\end{figure}

Incorporating these ideas using the consensus matrix $S_{6}$ from Section \ref{sec:dna} and an initial probability vector of all zeros except for a 1 in the second position gives us the custom cluster of $\{2,4,6,15,19,29,36\}$, a cluster with four other AML-B patients and two AML patients (although one of them, Patient 29, consistently clusters with the AML-B patients in our experience). These results are presented in table \ref{tbl:dnacc} along with the six nearest neighbors of Patient 2 using Euclidean distance and cosine measure. The SCA's custom cluster for Patient 2 features three patients not found in these nearest neighbor sets and suggests that physicians could learn a great deal by examining these hidden connections between Patient 2 and Patients 15, 29, and 36.

\begin{table}[htdp]
\caption[Custom Cluster for leukemia Patient 2]{Custom Cluster for leukemia Patient 2. This table shows the six other patients most similar to Patient 2. The patients are listed in similarity order, that is, the first one is the one most similar to Patient 2. The cluster returned by the SCA differs by three patients with both lists derived from two traditional distance measures.}
\begin{center}
\begin{tabular}{|c|c|} \hline
Method & Other Patients \\ \hline
SCA & $29,19,4,15,36,6$ \\
2-norm & $19,16,9,3,6,18$ \\
cosine & $16,19,9,3,18,4$ \\ \hline
\end{tabular}
\end{center}
\label{tbl:dnacc}
\end{table}

\section{Discussion}
These initial tests prove that the SCA can be an effective clustering tool. As with any new method, this initial promise raises multiple questions for further study, some of which are listed here.

\begin{itemize}
\item Use probabilistic analysis of initial probability vectors to see what we can learn about the number of possible solutions the SCA can return and whether there is any connection between $\sigma(P,n_{1})$ and the tendency of $P$ to produce multiple solutions.
\item Devise a fuzzy clustering for a data set based on the multiple results returned when using different initial probability vectors.
\item Investigate whether in situations where the stochastic clustering algorithm returns multiple answers, if building a consensus matrix from these results, and applying the SCA again will eventually yield a unique solution.
\item Examine whether the Sinkhorn-Knopp balancing step can be replaced by a simple scaling to make all row sums equal. Though we lose the results from Markov chain theory, perhaps they are unneeded since all we are looking for is $x_{t}^{T}$ values that are approximately equal. The work of Lin and Cohen mentioned in Section \ref{sec:singlemeasure} would seem to indicate that this is a possibility.
\item Continue the search for a single similarity measure whose values are distributed in a way that can be exploited by the stochastic clustering method.
\item Improve the bounds for values of $d_{i}$. Numerical results indicate that the upper bound found for Theorem \ref{th:sigma} can be greatly improved.  
\item Explore the structure of the spectrum of symmetric, irreducible, nearly uncoupled, doubly stochastic matrices. For this paper, we were only concerned with the eigenvalues near one, but from examining eigenvalues during the course of this research, there appears to be some structure to the spectrum, especially a large number of eigenvalues near zero.
\item Work to find a tighter bound on the numeric connection between $\lambda_{2}(P)$ and $\sigma(P,n_{1})$ that Theorems \ref{thm:hgeasy} and \ref{thm:hghard} establishes. 
\end{itemize}

\section{Acknowledgements} The authors would like to thank the referees for helpful comments that made this a much better paper. In particular two referees brought notable earlier results to our attention for which we are very grateful. Thanks also to Ilse Ipsen for her suggestions and the 2-norm bound in the proof of Theorem \ref{thm:hgeasy}.

\end{document}